\documentclass[12pt]{article}
\usepackage{amsfonts}
\begin{document}

\begin{center}
{\large\textbf{Fano's inequality is a mistake }}\end{center}
\begin{center}{\large{MARAT  GIZATULLIN }}\end{center}
\bigskip
\par\noindent The Department of Mathematics,\par\noindent
 Technical University Federico Santa Mar{\'{\i}}a,\par\noindent
Avenida Espa\~na, No. 1640, Casilla 110-V,\par\noindent
Valpara{\'{\i}}so, Chile \medskip
\par\noindent {\emph{e-mail}} \textsf{mgizatul@mat.utfsm.cl}
\bigskip
\begin{center}
{\large\textbf{Introduction }}\end{center}
\par
Let us consider the simplest Fano threefold , that is the
three-dimensional projective space $ \mathbb{P}^3$. For this
threefold, Fano's inequality looks as follows. \medskip \par For
any Cremona transformation
$$ f: \mathbb{P}^3--\rightarrow \mathbb{P}^3$$ defined by four
homogeneous polynomials of the same degree $d$ and without a
common non-constant factor,
$$x_{i}^{\prime}= f_{i}(x_0,x_1,x_2,x_3),\quad i=0,1,2,3,$$
either there exists a point $P\in \mathbb{P}^3$
 such that $${\textnormal{mult}}_P(f_{i})> d/2 $$ for every $i=0,1,2,3, $
 or there exists an irreducible curve $C\subset \mathbb{P}^3$
 such that $${\textnormal{mult}}_C(f_{i})> d/4 $$ for every $i$.\medskip \par
 One can remark that for the first time, these inequalities were
  indicated by Margherita Piazzola-Beloch  in [1]. She was a pupil of G. Castelnuovo,
 her paper presents the text of her thesis , G.  Castelnuovo
 was the adviser of the thesis. Thus all (including G. Fano ) the subsequent authors
 of the variants or generalizations of the Fano inequality are out of
 the historical responsibility for the mistake explained below.  \medskip
 \par The goal of my article is to show that these inequalities
 do not take place for a Cremona transformation of degree 13,
 that is I  write down the formulas for a Cremona
 transformation of degree 13 such that for the forms $f_0,...,f_4 $
 defining the transformation, for any point $P\in \mathbb{P}^3$
 and
 for any curve $C\subset \mathbb{P}^3 $ one can see that
$${\min}_i({\textnormal{mult}}_P(f_{i}))\leq 6 , \qquad {\min}_i
({\textnormal{mult}}_C(f_{i}))\leq 3. $$
\medskip
\begin{center}
{\large\textbf{The construction of the example}}\end{center}
\medskip
\par Let us consider the  homogeneous coordinates $x_0,x_1,x_2,x_3 $
for $ \mathbb{P}^3$ as the normalized coefficients of a binary
cubic form $F(T_0,T_1),$
$$F(T_0,T_1)=x_{0}T_{0}^3+3x_{1}T_{0}^2T_{1}+3x_{2}T_{0}T_{1}^2+x_{3}T_{1}^3.$$
Let $D=D(x_0,x_1,x_2,x_3)$ be the discriminant of the binary
cubic, $$D=
x_0^2x_3^2-3x_1^2x_2^2-6x_0x_1x_2x_3+4x_0x_2^3+4x_3x_1^3.$$ Let us
fix a parameter $t$ and consider four following forms of degree
13.
$$(f_t)_0=x_0D^3, $$
$$(f_t)_1=x_1D^3+tx_0^5D^2,$$
$$(f_t)_2=x_2D^3+2tx_1x_0^4D^2+t^2x_0^9D,$$
$$(f_t)_3=x_3D^3+3tx_2x_0^4D^2+3t^2x_1x_0^8D+t^3x_0^{13}.$$
These four forms define a one-parameter family of rational maps
$$ g_t :\mathbb{P}^3--\rightarrow \mathbb{P}^3.$$
 If $t=0,$ then all the four forms have a common factor,
the factor is $D^3,$ after the cancellation we see that $g_0$ is
the identity transformation. For our example we need non-zero
values of $t$. If $t$ is not zero, then it is clear that the four
forms are without a common non-constant factor. Further,
$$D\Big((f_t)_0,(f_t)_1,(f_t)_2,(f_t)_3\Big)=D(x_0,x_1,x_2,x_3)^{13},$$
this identity
actually expresses the invariant property (with respect to the
triangular transformation of variables $T_0,T_1 $ ) of the
discriminant. Using the latter identity, is is not hard to see
that
$$(f_{-t})_{i}\Big((f_t)_0,(f_t)_1,(f_t)_2,(f_t)_3\Big)=x_{i}D^{42}, $$
that is $$g_{(-t)}\circ g_t={\textnormal{the identity
transformation}}.$$ Thus $g_t$ is rationally invertible and is a
Cremona transformation. More generally,
$$g_{s}\circ g_t=g_{s+t},$$
and we get a one-parameter group of Cremona transformations. These
transformations induce biregular automorphisms of an affine open
subset of the projective space, the subset is the complement to
the discriminant quartic  surface $D=0. $ Indeed, the above
formula of the discriminant transformation proves it (moreover,
one can see below the exact calculation of the fundamental points
of such a transformation ) . \par The formulas for $g_t $ ( or
more general formulas for an infinite-dimensional family of
automorphisms of the complement to the discriminant surface) were
written down on page 8 of the Max-Planck-Institute preprint [2]
.\medskip \par Let us fix a nonzero value of the parameter $t,$
for example put $t=1,$ and consider the corresponding Cremona
transformation
$$x_0'=x_0D^3, $$
$$x_1'=x_1D^3+x_0^5D^2,$$
$$x_2'=x_2D^3+2x_1x_0^4D^2+x_0^9D,$$
$$x_3'=x_3D^3+3x_2x_0^4D^2+3x_1x_0^8D+x_0^{13}.$$
First of all, we will find the points $P$ where the multiplicities
of every of the right hand sides are positive ( that is the set of
all common zeros of these right hand sides, or , in other words,
the fundamental points of the transformation).\par The first right
hand side vanishes if either $x_0=0$,  or $D=0$, or simultaneously
$x_0=0, D=0.$ \par If $x_0=0,$ but $D\neq 0,$ then
 using other three formulas, one can see that for other  three coordinates of a
 fundamental point , the equalities $x_1=x_2=x_3=0$ take place,  that is
in this case we are out of the projective space. \par The case
$D=0,$ but $x_0\neq 0$ is also impossible for a fundamental point.
\par Thus , the set of fundamental points consists of the
solutions of the following system of equations
$$ x_0=0, \quad D=0,$$ or equivalently,
$$ x_0=0, \quad x_1^2(4x_3x_1-3x_2^2)=0.$$
We see that  the set of fundamental points is the union of two
curves, the first curve is line $L$, $$L \quad :\qquad  x_0=0,
\quad x_1=0,$$ the second curve is conic $C$,
$$ C \quad :\qquad x_0=0, \quad 4x_3x_1-3x_2^2=0.$$
The discriminant surface $D=0$ has double points disposed on the
twisted cubic $T$ having the following homogeneous
parameterization,
$$ x_0=t_0^3,\quad x_1=t_0^2t_1,\quad x_2=t_0t_1^2, \quad x_3=t_1^3.$$
More precisely, the singular locus of the discriminant surface is
$T,$ and  mult$_P(D)=2$ for every $ P\in T$. Therefore
 general points of the line $L$  and of the conic $C$  are the
 points of
multiplicity one on the discriminant,
$$ {\textnormal{mult}}_L(x_1D^3+x_0^5D^2)=3,$$
$$ {\textnormal{mult}}_C(x_1D^3+x_0^5D^2)=3.$$
More generally, if a point $P$ of the intersection of discriminant
surface with the plane $x_0=0 $
 is
located out of the twisted cubic $T$, then at least one of the
multiplicities $ {\textnormal{mult}}_P(x_iD^3)$ is equal to 3. The
last hope to get a point of higher multiplicity is to consider the
point of intersection of the twisted cubic $T$ with the union of
curves $ L$ and $C$. It is obvious that
$$ T\cap ( L\cup C)=T\cap  L\cap C =\{ Q\},\quad Q=(0:0:0:1),$$
but it clear that
$$
{\textnormal{mult}}_Q(x_3D^3+3x_2x_0^4D^2+3x_1x_0^8D+x_0^{13})$$
$$={\textnormal{mult}}_Q(D^3)= 6.$$

\end{document}